\documentclass[12pt]{extarticle}

\usepackage[margin=1.5in]{geometry}  
\usepackage{graphicx}              
\usepackage{amsmath}               
\usepackage{amsfonts}              
\usepackage{amsthm}                
\usepackage{framed}

\newtheorem{thm}{Theorem}[section]
\newtheorem{lem}[thm]{Lemma}

\newtheorem{cor}[thm]{Corollary}

\begin{document}

\nocite{*}

\title{\bf An Explicit Result for $|L(1+it,\chi)|$}

\author{\textsc{Adrian W. Dudek} \\ 
Mathematical Sciences Institute \\
The Australian National University \\ 
\texttt{adrian.dudek@anu.edu.au}}
\date{}

\maketitle

\begin{abstract}
\noindent We give an explicit upper bound for non-principal Dirichlet $L$-functions on the line $s=1+it$. This result can be applied to improve the error in the zero-counting formulae for these functions.
\end{abstract}

\section{Introduction}

Let $\chi$ be a non-principal Dirichlet character to the modulus $q \geq 1$ and let $L(s,\chi)$ denote its associated $L$-function. Our objective here is to establish some explicit bounds for $| L(1+it, \chi)|$; the main result is given in the following theorem.

\begin{thm} \label{one}
Let $q \geq 1$ be a non-principal modulus with associated Dirichlet character $\chi$. Then
$$| L(1+it, \chi) | < \frac{\varphi(q)}{q} \log t + \log q + \gamma$$
for all $t > 50$, where $\varphi(q)$ is Euler's totient function and $\gamma$ is Euler's constant.
\end{thm}

We prove this using a method due to Backlund \cite{backlund}. Of course, we would not expect the above theorem to hold for values of $t$ which are close to zero. To this end, we prove the following theorem using the method of partial summation.

\begin{thm} \label{two}
Let $q \geq 1$ be a non-principal modulus with associated Dirichlet character $\chi$. Then
$$| L(1+it, \chi) | < \log (t+14/5) + \log q + 1$$
for all $t > 0$.
\end{thm}
It is then a straightforward exercise to combine Theorems \ref{one} and \ref{two} to get a result which is as good as Theorem \ref{one} asymptotically but which holds for all $t>0$. Note that
$$\gamma < \log\bigg(e^{\gamma} + \frac{109}{2t}\bigg)$$
for all $t > 0$ and also that
$$1+ \log(t+14/5) < \log\bigg(e^{\gamma} + \frac{109}{2t}\bigg)$$
for all $0<t\leq 50$. This gives us the following result.

\begin{cor}
Let $q \geq 1$ be a non-principal modulus with associated Dirichlet character $\chi$. Then
$$| L(1+it, \chi) | < \log ( q (e^{\gamma} t+ 109/2) )$$
for all $t > 0$.
\end{cor}
These are not deep results, though the author could not find any explicit bounds in the literature. It should be noted that bounds of this type can be used to explicitly estimate $N(T,\chi)$, the number of zeroes of $L(s,\chi)$ in the critical strip with height at most $T$. One could use Theorems \ref{one} and \ref{two} to improve the recent work of Trudgian \cite{trudgianN} in this way. Bounds for $N(T,\chi)$ can then be used for explicit estimates involving the prime number theorem for arithmetic progressions; see the work of McCurley \cite{mccurley} for a demonstration of this.

We will prove Theorems \ref{one} and \ref{two} in the following two sections. We only need to consider $q \geq 2$, noting that the case $q=1$ corresponds to the Riemann zeta-function. Explicit bounds in this case have been given by Trudgian \cite{trudgianbound} and Ford \cite{ford}.

To construct such a bound for a Dirichlet $L$-function on the line $s=1+it$, we first consider the Hurwitz zeta-function  
$$\zeta(s,c) = \sum_{n=0}^{\infty} \frac{1}{(n+c)^s}$$
for $c \in (0,1]$ and $\Re(s) > 1$. We then have that
\begin{eqnarray} \label{connection}
L(s,\chi) & = & \sum_{n=1}^{\infty} \frac{\chi(n)}{n^s} \nonumber \\
& = & \sum_{1 \leq a \leq q} \sum_{n \equiv a \text{ mod } q} \frac{\chi(n)}{n^s} \nonumber \\
& = & \sum_{1 \leq a \leq q} \chi(a) \sum_{n \equiv a \text{ mod } q} \frac{1}{n^s} \nonumber \\
& = & \sum_{1 \leq a \leq q} \frac{\chi(a)}{q^s} \zeta\Big(s,\frac{a}{q}\Big) 
\end{eqnarray}
where the manipulations are justified by the absolute convergence of the Dirichlet series for $\text{Re}(s)>1$. The connection between $L(s,\chi)$ and $\zeta(s,c)$ is given hearty discussion in Chapter 12 of Apostol's \cite{apostol} introductory text, where he uses the above equation to define $L(s,\chi)$ for \textit{all} $s\in \mathbb{C}$ (where $\chi$ is, importantly, non-principal). It so follows from (\ref{connection}) that
\begin{equation} \label{trivial}
|L(s,\chi)| \leq \frac{1}{q}\sum_{\substack{1 \leq a \leq q \\ (a,q)=1}} \Big|\zeta\Big(s,\frac{a}{q}\Big)\Big|.
\end{equation}
Therefore, it remains to develop bounds for the Hurwitz zeta-function.

\section{The proof}

\subsection{Backlund's method}

Our technique for Theorem \ref{one} is that used by Backlund \cite{backlund} in bounding $\zeta(1+it)$. We start with the summation formula of Euler and Maclaurin given in Murty's \cite{murty} book of problems in analytic number theory.

\begin{lem} \label{EM}
Let $k$ be a nonnegative integer and $f(x)$ be $(k+1)$ times differentiable on the interval $[a,b]$. Then
\begin{eqnarray*}
\sum_{a < n \leq b} f(n) & = & \int_a^b f(t) dt +\sum_{r=0}^{k} \frac{(-1)^{r+1}}{(r+1)!} \Big\{ f^{(r)}(b) -f^{(r)}(a) \Big\} B_{r+1}\\
 & + & \frac{(-1)^k}{(k+1)!} \int_a^b B_{k+1}(x) f^{(k+1)}(x)dx,
\end{eqnarray*}
where $B_j(x)$ is the $j$th periodic Bernoulli polynomial and $B_j = B_j(0)$.

\end{lem}

We consider that for $\Re(s) > 1, t>0$ and some integer $N>0$ we have
$$\sum_{n > N} \frac{1}{(n+\frac{a}{q})^s} = \zeta\Big(s,\frac{a}{q}\Big) -\sum_{0 \leq n \leq N} \frac{1}{(n+\frac{a}{q})^s}. $$
We use Lemma \ref{EM} with $f(n) = (n+a/q)^{-s}$, $k=1$, $a = N$ and $b \rightarrow \infty$ to estimate the left side of the above formula, and so we have for the right side:
\begin{eqnarray*}
\zeta\Big(s,\frac{a}{q}\Big) -\sum_{0 \leq n \leq N-1} \frac{1}{(n+\frac{a}{q})^s} & = & \frac{(N+\frac{a}{q})^{1-s}}{s-1} - \frac{1}{2 (N+\frac{a}{q})^{s}} + \frac{s}{12 (N+\frac{a}{q})^{s+1}} - \\
& - & \frac{s(s+1)}{2} \int_N^{\infty} \frac{ \{x\}^2 - \{x\} + 1/6}{(x+\frac{a}{q})^{s+2}}dx.
\end{eqnarray*}
Note that as the right side of the above equation remains valid for $\Re(s) > -1$, it must remain valid for $s = 1+it$ where $t>0$. Thus, setting $N = [(t/m)-(a/q)]$, where $m$ is to be chosen later, we obtain the estimate
$$|\zeta(1+it,a/q)| - \sum_{0 \leq n \leq N-1} \frac{1}{n+a/q} \leq \frac{1}{t} + \frac{m}{2(t-m)} + \frac{m^2 (1+t)(4+t)}{24 (t-m)^2}.$$
The sum can be estimated using the classic bound
\begin{equation} \label{classic}
\sum_{1 \leq n \leq t} \frac{1}{n} \leq \log t + \gamma + \frac{1}{t}
\end{equation} 
where $\gamma$ is the Euler--Mascheroni constant. We get
\begin{eqnarray*}
\sum_{0 \leq n \leq N-1} \frac{1}{n+a/q} & \leq & \frac{q}{a} + \sum_{1 \leq n \leq N} \frac{1}{n} - \frac{1}{N} \\
& \leq & \log N + \gamma + \frac{q}{a} \\
& \leq & \log t - \log m + \gamma + \frac{q}{a}.
\end{eqnarray*}

Putting everything together we have
$$|\zeta(1+it,a/q)| - \log t - \frac{q}{a} \leq - \log m  + \gamma + \frac{1}{t} + \frac{m}{2(t-m)} + \frac{m^2 (1+t)(4+t)}{24 (t-m)^2}.$$
If we now choose $m=3$ it is easy to verify that the right side of the above inequality is negative for all $t \geq 50$. Therefore, we have that
$$|\zeta(1+it,a/q)| < \log t + \frac{q}{a}.$$
We note that we could have replaced $a/q$ with some real number $c\in(0,1]$ in the above working; this gives us the following result for the Hurwitz zeta-function.

\begin{lem}
Let $c \in (0,1]$. Then
$$|\zeta(1+it,c)| < \log t + \frac{1}{c}.$$
\end{lem}

From (\ref{trivial}) we get
\begin{eqnarray*} 
|L(1+it,\chi)| & \leq & \frac{1}{q} \sum_{\substack{1 \leq a \leq q \\ (a,q)=1}} \bigg( \log t + \frac{q}{a}\bigg) \\
& \leq & \frac{\varphi(q)}{q} \log t +  \sum_{\substack{1 \leq a \leq q \\ (a,q)=1}} \frac{1}{a}.
\end{eqnarray*}
It remains to use (\ref{classic}) to estimate the sum as follows.
\begin{eqnarray*}
\sum_{\substack{1 \leq a \leq q \\ (a,q)=1}} \frac{1}{a} & \leq & -\frac{1}{q} + \sum_{1 \leq a \leq q} \frac{1}{a}  \\
& \leq & \log q + \gamma. 
\end{eqnarray*}
This proves Theorem \ref{one}. The reader should note that there is a more careful way in which one could estimate the above sum (see, for example, Exercise 1.5.1 of Murty \cite{murty}). This would introduce a factor of $\varphi(q)/q$ to the front of the logarithm term along with additional error terms. We omit this for the sake of brevity. 

\subsection{By partial summation}

As before, let $\chi$ be a non-principal character of modulus $q$. Using Abel's method of partial summation we have
$$\sum_{n=1}^N \frac{\chi(n)}{n^{1+it}} = \Big( \sum_{n=1}^N \chi(n) \Big) N^{-1-it} + (1+it) \int_1^N \Big( \sum_{n \leq t} \chi(n) \Big) \frac{dt}{t^{2+it}}.$$
For ease of notation, we write 
$$A(t) = \sum_{n \leq t} \chi(n)$$
to get
\begin{equation} \label{bla}
\sum_{n=1}^N \frac{\chi(n)}{n^s} = A(N) N^{-s} + s \int_1^N \frac{A(t)}{t^{s+1}}dt.
\end{equation}
Note that we can manipulate $L(s,\chi)$ as follows
\begin{eqnarray*}
L(s,\chi) & = & \sum_{n=1}^{\infty} \frac{\chi(n)}{n^s} = \sum_{n=1}^{\infty} \Big( A(n) - A(n-1) \Big) \frac{1}{n^s} \\
& = & \sum_{n=1}^{\infty} A(n) \Big( \frac{1}{n^s} - \frac{1}{(n+1)^s} \Big)
\end{eqnarray*}
to get
\begin{equation} \label{blu}
L(s,\chi) =  s \int_1^{\infty} \frac{A(t)}{t^{s+1}}dt.
\end{equation}
Subtracting (\ref{bla}) from (\ref{blu}) we have
$$L(s,\chi) - \sum_{n=1}^N \frac{\chi(n)}{n^s} = s \int_N^{\infty} \frac{A(t)}{t^{s+1}}dt - A(N) N^{-s}.$$
Clearly, $|A(t)| < q/2$, as $A(N)$ is equal to zero whenever $N$ is a multiple of $q$. The absolute convergence of the right side for $\Re(s) > 0$ allows us to set $s=1+it$ to get
$$|L(1+it,\chi)| \leq  \sum_{n=1}^N \frac{1}{n} + \frac{q(1+t)}{2} \int_N^{\infty} \frac{1}{t^{2}}dt + \frac{A(N)}{N}.$$
If we set $N = [q(t+b)/m]$ (where $m$ and $b$ are to be chosen later) and use (\ref{classic}) we have
$$|L(1+it,\chi)| - \log(q(t+b))-1 \leq - \log m + (\gamma-1)  + \frac{m(2+q+qt)}{2 q(t+b)-2m}.$$
If we choose $m=2$ and $b=14/5$, it is straightforward to verify that the right hand side is negative for all $q \geq 2$ and $t >0$. This completes our proof.

\clearpage

\bibliographystyle{plain}

\bibliography{biblio}

\end{document}